\newcommand{\rem}[1]{}
\documentclass[10pt,oneside,reqno]{amsart}
\usepackage{amsfonts,amssymb,amsmath,amsthm,mathrsfs}
\usepackage{url}
\usepackage{graphicx}
\usepackage{epstopdf}
\theoremstyle{definition}

\begin{document}
\author[L.~G.~Molinari]{Luca~Guido~Molinari}
\address{Physics Department ``Aldo Pontremoli'', Universit\`a degli Studi di Milano and I.N.F.N. sez. Milano,  
Via Celoria 16, 20133 Milano, Italy. ORCiD (L.G.Molinari): 0000-0002-5023-787X}
\email{luca.molinari@unimi.it}
\title[Trigonometric sums]{A note on trigonometric approximations \\of Bessel functions of the first kind,\\
and trigonometric power sums}
\begin{abstract} 
I reconsider the approximation of Bessel functions with finite sums of trigonometric functions, in the light 
of recent evaluations of Neumann-Bessel series with trigonometric coefficients. A proper choice of the angle allows for
an efficient choice of the trigonometric sum.
Based on these series, I also obtain straightforward non-standard 
evaluations of new parametric sums with powers of cosine and sine functions. 
\end{abstract}
\date{21 May 2022}
\maketitle
\section{Introduction}
Bessel functions are among the most useful and studied special functions. Analytic expansions exist for different regimes \cite{NIST}, and numerical algorithms for their 
precise evaluation \cite{Matviyenko93}\cite{Schwartz12}\cite{Bremer17}\cite{Karatsuba19}. 
Their simplest approximations are polynomials \cite{Abramowitz}\cite{Gross95}\cite{Li06} 
and finite trigonometric sums, that can be advantageous in applications \cite{Barsan07}. \\
Let's consider $J_0$. Several trigonometric sums appeared in the decades, sometimes being rediscovered. These very simple ones
\begin{align}
& J_0(x) \simeq \tfrac{1}{4}[1+\cos x + 2\cos (\tfrac{\sqrt 2}{2}x)] \label{APP1}\\
& J_0(x)\simeq \tfrac{1}{6}\left [1+ \cos x + 2\cos (\tfrac{1}{2}x) + 2\cos (\tfrac{\sqrt 3}{2}x)  \right ] \label{APP2}
\end{align}
have errors $\epsilon = J_0-J_0^{\rm approx}$ with power series (the marvel of Mathematica)
$$ \epsilon (x) = -\tfrac{x^8}{2^8 \cdot 20560} (1-\tfrac{x^2}{36}+ ... ), \qquad   \epsilon(x) = -\tfrac{x^{12}}{2^{12}\cdot 239500800}
(1-\tfrac{x^2}{52}+ ...) $$
In practice, an error less than $0.001$ is achieved for $x \le 3$ or $x\le 5.9$. 
These approximations were obtained by Fettis with the Poisson formula \cite{Fettis}. Rehwald \cite{Rehwald67} and later  
Waldron \cite{Waldron81}, Blachman and Mousavinezhad \cite{Blachman86} and \cite{Abuelmaatti} used the strategy of truncating to the first term Neumann-Bessel series like 
$$J_0(x)+2J_8(x) + 2J_{16}(x) + ... = \tfrac{1}{4} \cos [1+\cos x + 2\cos (\tfrac{\sqrt 2}{2}x)] $$ 
that can be obtained from the Bessel generating function. The examples correspond to $n=4,6$ of eq.19 in \cite{AlJarrah02}:
\begin{align}
J_0(x)+2{\sum}_{k=1}^{\infty}  (-)^{kn} J_{2kn}(x) =\frac{1}{n} {\sum}_{\ell=0}^{n-1} \cos ( x\cos \tfrac{\pi}{n}\ell ) \label{SUM}
\end{align}
and the errors reflect the behaviour $J_{2n}(x)\approx (x/2)^{2n}$ of the first neglected term, but with much larger denominators. \\
The truncation yields $J_0$ as a sum of cosines that corresponds to the evaluation of the Bessel integral $J_0(x)= \int_0^\pi \frac{d\theta}{\pi} \cos (x\cos \theta )$ with the trapezoidal rule with $n$ nodes 
\cite{Stroud69}\cite{Baratella81}\cite{Trefethen14}. Increasing $n$ increases accuracy: $n=15$ is a formula by Fettis \cite{Fettis} with 8 cosines (instead of 15, symmetries of the roots of unity reduce the number of terms):
\begin{align}
J_0(x)\simeq \tfrac{1}{15}  \cos x + \tfrac{2}{15} {\sum}_{k=1}^7 \cos (x \cos \tfrac{k\pi}{15})  \label{APP3}
\end{align}
The error now is order $x^{30}\times 10^{-42}$ and less than $10^{-6}$ for $x<15$.

In this note I reconsider the approximations for $J_0$  in the light of new Neumann-Bessel trigonometric 
series in ref.\cite{Molinari21}. They extend the series \eqref{SUM} by including an angular parameter, that is chosen to kill the 
term with $J_{2n}$, so that the truncation involves the next-to-next term $J_{4n}$ of the series.\\
The same strategy on appropriate series is then used for Bessel functions $J_n$ of low order, that are discussed in section 3. \\
In section 4, I show that the series give in very simple way some parametric sums of powers of cosines and sines.  Some are
in the recent literature (Jelitto \cite{Jelitto22}, 2022) while the following ones, to my knowledge, are new:
$$ \sum_{\ell=0}^n \sin^p (\tfrac{\theta+2\pi\ell}{n}) \begin{array}{c} \sin\\ \cos \\ \end{array} (q \tfrac{\theta+2\pi\ell}{n}) \qquad (p,q=0,1,...).$$

\section{The Bessel function $J_0$} 
Consider the Neumann trigonometric series eq.11 in \cite{Molinari21}: 
\begin{align}
J_0(x)+2\sum_{k=1}^{\infty}  (-)^{kn} J_{2kn}(x)  \cos(2kn\theta)
=\frac{1}{n} \sum_{\ell=0}^{n-1} \cos [ x\cos (\theta +\tfrac{\pi}{n}\ell )]  \label{sum1}
\end{align}
The approximations \eqref{APP1}, \eqref{APP2} and \eqref{APP3} are obtained with $\theta = 0$, $n=4$, $6$, $15$, and neglecting functions $J_8$, $J_{12}$, $J_{30}$ and higher orders.  However they are not optimal.  The advantage of eq.\eqref{sum1} 
is the possibility to choose the angle $\theta=\pi/{4n}$ to kill all terms $J_{2n}$, $J_{6n}$, etc. Then:
\begin{align}
J_0(x) - 2 J_{4n}(x) + 2J_{8n} (x) - ... =\frac{1}{n} \sum_{\ell=0}^{n-1} \cos ( z\cos \tfrac{1+4\ell}{4n}\pi  )  \label{1p2}
\end{align}
An expansion for $J_0$ results, again, by neglecting the other terms.\\ 
Some examples:\\
$\bullet $ $n=2$. It is 
$J_0(x)=\tfrac{1}{2}[\cos (x\cos \tfrac{\pi}{8}) +\cos (x\sin\tfrac{\pi}{8})] + \epsilon_2(x)$. If we neglect the error,
the first zero occurs at  $\pi \sqrt{2-\sqrt 2}=2.4045$ ($j_{0,1}=2.4048$). \\
$\bullet $ $n=3.$ The approximation has three cosines:
\begin{align}
J_0(x) =& 
\tfrac{1}{3} [\cos (x\tfrac{1}{\sqrt 2}) + \cos (x \tfrac{\sqrt 3-1}{2\sqrt 2}) + \cos (x \tfrac{\sqrt 3 +1}{2\sqrt 2})] +
\epsilon_3 (x)\label{J03}\\
 \epsilon_3 (x)=& \tfrac{x^{12}}{2^{12}\cdot 239500800}\left[ 1 - \tfrac{x^2}{52} + \tfrac{x^4}{52\cdot 112}-\tfrac{x^6}{52\cdot 112\cdot 180} + \ldots \right ]. \nonumber
 \end{align}
Remarkably, the first powers of the error are opposite of those for the expansion eq.\eqref{APP2}, that would involve 6 terms if
not for the degeneracy of the roots of unity. 
The half-sum of \eqref{APP2} and \eqref{J03},
\begin{align}
J_0(x)\simeq  \tfrac{1}{12}\big [1+ \cos x + 2\cos (\tfrac{1}{2}x) + 2\cos (\tfrac{\sqrt 3}{2}x) + 2\cos (x\tfrac{1}{\sqrt 2})  \nonumber \\
+2\cos (x \tfrac{\sqrt 3-1}{2\sqrt 2}) + 2\cos (x \tfrac{\sqrt 3 +1}{2\sqrt 2})\big ] \label{eps24}
\end{align} 
has error 
$\epsilon (x) = -\tfrac{x^{24}}{5.2047} \times 10^{-30} [ 1- \tfrac{x^2}{100}+ \tfrac{x^4}{20800} - \dots ] $.\\
%
%
%
$\bullet $ $n=6$ gives a precision similar to the sum \eqref{eps24}:
\begin{align}
J_0 (z) = &
\tfrac{1}{6}\big [\cos (x \cos \tfrac{\pi}{24}) +\cos ( x\cos \tfrac{3\pi }{24})+\cos (x\cos \tfrac{5\pi}{24}  ) \nonumber\\
&+\cos (x \sin \tfrac{\pi}{24}) +\cos ( x\sin\tfrac{3\pi }{24})+\cos (x\sin \tfrac{5\pi}{24}  ) \big ] +\epsilon_6 (x). \label{J0n6}
\end{align}
The error has power expansion $ \epsilon_6 (x) = \frac{x^{24}}{5.2047} \times 10^{-30} [1-\frac{x^2}{100} + \ldots ]$.\\
$\bullet $ $n=8$ is a sum of 8 cosines and compares with the formula \eqref{APP3} by Fettis.
The two approximations are different but with the same number of terms (because $\theta=0$ produces degenerate terms) and similar
precision.

\begin{figure}[t]
\begin{center}
\includegraphics[width=7cm]{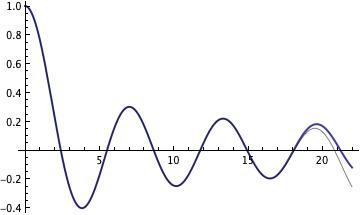}
\hspace{1cm} 
\caption{\label{figure1:fig} The Bessel function $J_0$ (thick) and the trigonometric expansion \eqref{J0n6}. The difference increases with $x$;
it is less than $10^{-9}$ for $x<8$ and $10^{-3}$ for $x<15$.}
\end{center}
\end{figure}

\section{Bessel functions $J_n$.}
\noindent
$\bullet $ $\mathbf{J_1}$ is evaluated via $J_1=-J_0'$. Eq. \eqref{eps24} gives:
\begin{align}
J_1(x) \simeq  \tfrac{1}{12}\big [\sin x + \sin (\tfrac{1}{2}x) + \sqrt 3 \sin (\tfrac{\sqrt 3}{2}x) + 
\sqrt 2 \sin(x\tfrac{1}{\sqrt 2})  \nonumber \\
+(\tfrac{\sqrt 3-1}{\sqrt 2}) \sin (x \tfrac{\sqrt 3-1}{2\sqrt 2}) 
+(\tfrac{\sqrt 3+1}{\sqrt 2})\sin (x \tfrac{\sqrt 3 +1}{2\sqrt 2})\big ] 
\end{align} 
with error  $\epsilon(x) \simeq (x/20)^{23} \times 3.87 \times [1-\tfrac{13}{1200}x^2 + \dots ].$\\
$\bullet$ $\mathbf{J_2, J_4}$ can be evaluated with the following identity
(the real part of eq.(5) in \cite{Molinari21}):
\begin{align}
J_p(x) +{\sum}_{k=1}^\infty [J_{kn+p} (x)+(-1)^{kn+p}J_{kn-p}(x)] \cos(kn\theta)  \nonumber\\
=\tfrac{1}{n}{\sum}_{\ell=0}^{n-1} \cos[x \sin(\theta+\tfrac{2\pi\ell}{n} ) + p (\theta+\tfrac{2\pi\ell}{n} )]\label{sum2}
\end{align}
Because of the term $J_{n-p}$, we take $2p<n$. With $y=\frac{\pi}{2n}$:
\begin{align*}
J_p(x)  - (-1)^p J_{2n-p}(x) +\ldots =\tfrac{1}{n}{\sum}_{\ell=0}^{n-1} \cos[x \sin(\tfrac{1+4\ell}{2n}\pi ) + p\tfrac{1+4\ell}{2n} \pi ]
\end{align*}
If only $J_p$ is kept, the approximation depends on the parity of $p$:
\begin{align}
J_p(x)  \simeq
\begin{cases}\quad 
 \cos ( p\tfrac{1+4\ell}{2n} \pi ) \times \tfrac{1}{n}{\sum}_{\ell=0}^{n-1} \cos[x \sin(\tfrac{1+4\ell}{2n}\pi )]  & \text{$p$ even} \\
\, -\sin ( p\tfrac{1+4\ell}{2n} \pi ) \times \tfrac{1}{n}{\sum}_{\ell=0}^{n-1} \sin[x \sin(\tfrac{1+4\ell}{2n}\pi )] & \text{$p$ odd} 
\end{cases}
\end{align}
$p=2$, $n=6$,  give the short formula
\begin{align}
J_2(x) \simeq 
\tfrac{1}{2\sqrt 3}\big [ \cos(x \sin \tfrac{\pi}{12})  - \cos( x \cos \tfrac{\pi}{12})  \big ]
\end{align}
with error $\epsilon (x)  = 2.69114 \times (x/10)^{-10}[1-\frac{x^2}{44}+ ...]$. 
The first zero is evaluated $\frac{2}{3}\pi \sqrt 6\simeq  5.1302$
($j_{2,1}=5.13562$). A better approximation is $n=8$, $y=\frac{\pi}{16}$:
\begin{align}
J_2(x) \simeq & \tfrac{1}{4} \cos(\tfrac{\pi}{8}) [\cos (x \sin \tfrac{\pi}{16} ) - \cos (x \cos \tfrac{\pi}{16}) ] \\
&+ \tfrac{1}{4}\sin (\tfrac{\pi}{8}) [ \cos (x \cos \tfrac{5\pi}{16} )   - \cos (x \sin \tfrac{5\pi}{16} )] \nonumber 
\end{align}
with error $\epsilon (x) = 7.00119 \times 10^{-16} x^{14} [1-\frac{x^2}{60}+...]$; $\epsilon(5)=3\times 10^{-6}$, $\epsilon (8) = 0.0010$.

For $J_4$ we select $p=4$, $n=8$, $\theta = \frac{\pi}{16}$. Now the lowest neglected term is $J_{12}$:
\begin{align}
J_4(x) \simeq 
\tfrac{\sqrt 2 }{8}\big [ & \cos(x \sin \tfrac{\pi}{16}) + \cos(x \cos \tfrac{\pi}{16})  
- \cos(x \sin \tfrac{5\pi}{16} ) - \cos(x \cos \tfrac{5\pi}{16})  \big ] 
\end{align}
The error is less that $10^{-3}$ at $x<6.3$.\\
$\bullet $ $\mathbf{J_3, J_5}$. A useful sum for odd-order Bessel functions is eq.(17) in \cite{Molinari21}:
\begin{align}
\sum_{k=0}^{\infty} (-)^{n+k} J_{(2n+1)(2k+1)}(x)  \cos[(2k+1)\theta ]
=\sum_{\ell=0}^{2n} \frac{\sin [x\cos (\tfrac{\theta+2\pi\ell}{2n+1} )]}{2(2n+1)} \label{17}
\end{align}
The angle $\theta= \frac{\pi}{6}$ cancels $J_{6n+3}$, $J_{14n+7}$ etc. and gives the approximation
\begin{align}
 J_{2n+1}(x)   \simeq \frac{(-1)^n}{\sqrt 3} \sum_{\ell=0}^{2n} \frac{\sin [x\cos \tfrac{1+12\ell}{12n+6}\pi ]}{2n+1} 
\end{align}
that neglects $J_{10n+5}$ etc. 
With  $n=1$ and $n=2$ we obtain:
\begin{align}
J_3(x)\simeq &- \tfrac{1}{3\sqrt 3} [ \sin(x  \cos \tfrac{\pi}{18}) - \sin ( x  \sin \tfrac{2 \pi}{9}) - \sin ( x  \sin  \tfrac{\pi}{9})] \label{J3}\\
J_5(x)\simeq & \tfrac{1}{5\sqrt 3} \big [\sin (x \cos \tfrac{\pi}{30}) +\sin (x \sin \tfrac{\pi}{15}) -
\sin (\tfrac{\sqrt 3}{2}x) \nonumber \\
&\qquad -\sin (x \sin \tfrac{4\pi}{15}) +
\sin (x \cos \tfrac{2\pi}{15})\big ].
 \end{align}
The expansion for $J_3$ has error $\epsilon =  2.33373\times 10^{-17} x^{15} [1-\tfrac{x^2}{64}+ ...]$.
The second one has error $ \epsilon = 1.92134 \, x^{25} \times 10^{-33}\times [1-\frac{x^2}{104} + ... ]$.

\begin{figure}
\begin{center}
\includegraphics[width=7cm]{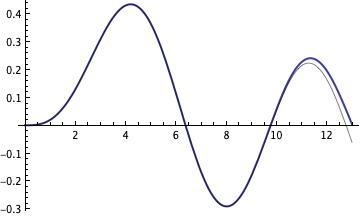}
\hspace{1cm} 
\caption{\label{figure1:fig} The Bessel function $J_3$ (thick) and the approximation \eqref{J3}. The difference is $\epsilon(6)=6\times 10^{-6}$,
$\epsilon(8)=.0003$, $\epsilon(10)=.0045$.}
\end{center}
\end{figure}

\section{Trigonometric identities}
The Neumann-Bessel series here used provide sums of powers of sines and cosines.
They arise by expanding in powers of $x$ the Bessel functions in the series, 
$$ J_n(x) =  \sum_{k=0}^\infty (-1)^k \frac{(x/2)^{n+2k}}{k! (k+n)!} $$
and the trigonometric functions in the sum of the series. \\

\noindent
$\bullet$ Consider the series eq.\eqref{sum1}. 
At threshold powers $x^{2n}$, $x^{4n}$ etc. new Bessel functions $ (-)^n 2 J_{2n}\cos(2n\theta)$, $2 J_{4n}\cos(4n\theta)$ etc. enter 
a term in the sum of cosines.
\begin{align}
& \tfrac{1}{n} {\sum}_{\ell=0}^{n-1} \left [\cos \tfrac{\theta+\ell\pi}{n}\right ]^{2k} =\\
& =\begin{cases}
\frac{1}{4^k}\binom{2k}{k} & 0\le k<n  \nonumber\\
\frac{1}{4^k} \left [\binom{2k}{k} +  2\binom{2k}{k-n} \cos(2\theta) \right ] & n\le k<2n \nonumber \\
\frac{1}{4^k} \left [\binom{2k}{k} +  2\binom{2k}{k-n} \cos(2\theta) +2 \binom{2k}{k-2n} \cos(4\theta) \right ] & 2n\le k< 3n \nonumber \\
... & ...
\end{cases}
\end{align}
By replacing $\theta $ with $\theta+n\frac{\pi}{n}$ we obtain:
\begin{align}
& \tfrac{1}{n} {\sum}_{\ell=0}^{n-1} \left [\sin \tfrac{\theta+\ell\pi}{n}\right ]^{2k} =\\
& =\begin{cases}
\frac{1}{4^k}\binom{2k}{k} & 0\le k<n  \nonumber\\
\frac{1}{4^k} \left [\binom{2k}{k} +  (-)^n 2\binom{2k}{k-n} \cos(2\theta) \right ] & n\le k<2n \nonumber \\
\frac{1}{4^k} \left [\binom{2k}{k} +  (-)^n 2\binom{2k}{k-n} \cos(2\theta) +2 \binom{2k}{k-2n} \cos(4\theta) \right ] & 2n\le k< 3n \nonumber \\
... & ...
\end{cases}
\end{align}
Examples: \\
$\frac{1}{9} {\sum}_{\ell=0}^8 \left [\sin  \frac{\theta+\ell\pi}{9} \right ]^{20}=\frac{1}{4^{10}} \left [\binom{20}{10} -  2\binom{20}{1} \cos(2\theta) \right ]. $\\
$\frac{1}{n} {\sum}_{\ell=0}^{n-1} \left [\cos ( \frac{1+6\ell}{6n} \pi )\right ]^{2n}=\frac{1}{4^n} \left [\binom{2n}{n} + 1\right ], $ \quad
$\frac{1}{n} {\sum}_{\ell=0}^{n-1} \left [\cos ( \frac{1+4\ell}{4n} \pi )\right ]^{4n}=\frac{1}{16^n} \left [\binom{4n}{2n} -  2 \right ]. $\\

\noindent
For  $\theta =0$ and $\theta=\frac{\pi}{2}$ these identities are eqs. 4.4.2 in \cite{Prudnikov}, 2.1 and 2.2 (together with several other non-parametric sums) in \cite{DaFonseca17}. The series had also been studied in \cite{Merca12}. 
Parametric averages on the full circle were recently evaluated by Jelitto \cite{Jelitto22}, with a  different method.\\ 

\noindent
$\bullet$ With the Neumann series \eqref{17} we obtain:
\begin{align}
& \tfrac{1}{2n+1} {\sum}_{\ell=0}^{2n} \left[\cos \tfrac{\theta+2\pi\ell}{2n+1} \right ]^{2k+1}=\\
&= \begin{cases}
 0 & \qquad\qquad  1\le 2k+1<2n+1 \nonumber\\
 \frac{1}{4^k} \binom{2k+1}{k-n} \cos \theta & \quad 2n+1\le 2k+1< 3(2n+1) \nonumber \\
 \frac{1}{4^k} \left [ \binom{2k+1}{k-n} \cos \theta +\binom{2k+1}{k-3n-1}\cos(3\theta )\right ]  & 3(2n+1)\le 2k+1< 5(2n+1) \nonumber \\
... & ... \nonumber
 \end{cases}
\end{align}
The sums of even powers of cosines are obtained from the series eq.16 in \cite{Molinari21}:
\begin{align}
J_0(x)+2\sum_{k=1}^{\infty} (-)^k J_{(4n+2)k}(x)  \cos (2k\theta )
= \sum_{\ell=0}^{2n} \frac{\cos[x\cos \tfrac{\theta+2\pi\ell}{2n+1}]}{2n+1}
\end{align}
\begin{align}
&\tfrac{1}{2n+1}{\sum}_{\ell=0}^{2n}\left [ \cos  \tfrac{\theta+2\pi\ell}{2n+1} \right ]^{2k}= \\
&=\begin{cases}
\frac{1}{4^k}\binom{2k}{k} & \qquad \; 0\le k<2n+1\\ 
\frac{1}{4^{k}}\left [\binom{2k}{k} +2 \binom{2k}{k-2n-1}\cos(2\theta)\right ] & 2n+1\le k<4n+2\\
\frac{1}{4^{k}}\left [\binom{2k}{k} +2 \binom{2k}{k-2n-1}\cos(2\theta) +2\binom{2k}{k-4n-2}\cos(4\theta) \right ] & 4n+2\le k<6n+3\\
... & ...
\end{cases} \nonumber
\end{align}
Example:
$(\cos \tfrac{\theta}{3})^{12} + (\cos \tfrac{\theta+\pi}{3})^{12}+ (\cos \tfrac{\theta+2\pi}{3})^{12} = \tfrac{3}{4^6} [\binom{12}{6}+2\binom{12}{3} \cos(2\theta) + 2\cos(4\theta)].$\\
The formulae with sines are obtained by shifting the parameter $\theta$.\\

\noindent
$\bullet$ Now let's consider the sum eq.\eqref{sum2} with $p<n-p$. 
The equations are new and are easier to state if we distinguish the parity of $n$ and of $p$.\\

\noindent
{\bf Case  $\mathbf {n=2m}$ and $\mathbf {p=2q}$}. Eq.\eqref{sum2} now is:
\begin{align*}
&\tfrac{1}{2m}{\sum}_{\ell=0}^{2m-1} \cos[x \sin (\tfrac{\theta+\pi\ell}{m} )] \cos[2q \tfrac{\theta+\pi\ell}{m}]= J_{2q}(x)+ \\
& + [J_{2m-2q}(x) + J_{2m+2q}(x)] \cos(2\theta ) +[J_{4m-2q}(x)+J_{4m+2q}(x)]\cos(4\theta)+...
\end{align*}
Separation of even and odd parts in $x$, and expansion in $x$ give:
$$\tfrac{1}{2m}{\sum}_{\ell=0}^{2m-1} [\sin \tfrac{\theta+\pi\ell}{m}]^{2k+1} \sin (2q \tfrac{\theta+\pi\ell}{m}) =0,\qquad \forall k$$ 
This result is obvious as the sum from $0$ to $m-1$ is opposite of the rest of the sum. The symmetry is used also
in the other result:
\begin{align}
&\tfrac{1}{m}{\sum}_{\ell=0}^{m-1} [\sin \tfrac{\theta+\pi\ell}{m}]^{2k} \cos (2q \tfrac{\theta+\pi\ell}{m}) =\\
&= \tfrac{(-)^q}{4^k} \begin{cases} 
0 & k<q \nonumber\\
\binom{2k}{k-q} & q\le k< m-q \nonumber \\
\binom{2k}{k-q} + (-1)^m \binom{2k}{k-m+q}\cos(2\theta) & m-q\le k< m+q \nonumber\\
\binom{2k}{k-q} + (-1)^m \left[ \binom{2k}{k-m+q}+ \binom{2k}{k-m-q} \right ] \cos(2\theta) & m+q\le k< 2m-q \nonumber\\
... & ... \nonumber
\end{cases} 
\end{align}

\noindent
{\bf Case $\mathbf{n=2m}$, $\mathbf{p=2q+1}$}. Eq.\eqref{sum2} becomes:
\begin{align*}
&-\tfrac{1}{2m}{\sum}_{\ell=0}^{2m-1} \sin[x \sin (\tfrac{\theta+\pi\ell}{m} )]\sin[ (2q+1) \tfrac{\theta+\pi\ell}{m} ] =J_{2q+1}(x)+\\
& + [-J_{2m-2q-1}(x) + J_{2m+2q+1}(x)] \cos(2\theta )
+[J_{4m-2q-1}(x)+J_{4m+2q+1}(x)]\cos(4\theta)+...
\end{align*}
The non trivial result is:
\begin{align}
&\tfrac{1}{m}{\sum}_{\ell=0}^{m-1} [\sin \tfrac{\theta+\pi\ell}{m}]^{2k+1} \sin [(2q+1) \tfrac{\theta+\pi\ell}{m}] =\\
&=\tfrac{(-)^q}{2^{2k+1}} 
\begin{cases} 0 & k<q \nonumber\\
\binom{2k+1}{k-q} & q\le k< m-q-1 \nonumber \\
 \binom{2k+1}{k-q} + (-1)^m \binom{2k+1}{k+m-q}\cos(2\theta) & m-q-1\le k< m+q \nonumber\\
 \binom{2k+1}{k-q} + (-1)^m \left [\binom{2k+1}{k+m-q}+ \binom{2k+1}{k-m-q} \right ] \cos(2\theta) & m+q\le k< 2m-q -1\nonumber\\
... & ... \nonumber
\end{cases} 
\end{align}
Example: $ \frac{1}{5}{\sum}_{\ell=0}^4 \sin^{13}(\tfrac{\pi\ell}{5}) \sin (\tfrac{3\pi\ell}{5}) = -
\frac{1}{2^{13}} \left [ \binom{13}{5} - \binom{13}{10} - \binom{13}{0} \right ] = -\frac{125}{1024}$.\\

\noindent
{\bf Case $\mathbf{n=2m+1}$ and $\mathbf{p=2q}$}:
\begin{align*}
\tfrac{1}{2m+1}{\sum}_{\ell=0}^{2m} \cos [x \sin (\tfrac{\theta+2\pi\ell}{2m+1} )] \cos (2q \tfrac{\theta+2\pi\ell}{2m+1}) &=J_{2q}(x)+
[J_{4m+2-2q}(x)+J_{4m+2+2q}(x)]\cos(2\theta)+...\\
\tfrac{1}{2m+1}{\sum}_{\ell=0}^{2m} \sin[x \sin (\tfrac{\theta+2\pi\ell}{2m+1} )] \sin(2q \tfrac{\theta+2\pi\ell}{2m+1}) &=
  [J_{2m+1-2q}(x) - J_{2m+1+2q}(x)] \cos \theta + ...
\end{align*}
\begin{align}
&\tfrac{1}{2m+1}{\sum}_{\ell=0}^{2m} [\sin \tfrac{\theta+2\pi\ell}{2m+1}]^{2k} \cos (2q \tfrac{\theta+2\pi\ell}{2m+1}) =\\
&= \tfrac{(-)^q}{4^k} \begin{cases} 
  0  &   k< q \\
 \binom{2k}{k-q} & q\le k< 2m+1-q \nonumber \\
 \binom{2k}{k-q} - \binom{2k}{k+q-1-2m} \cos (2\theta) & 2m+1-q \le k< 2m+1+q \nonumber \\
 \binom{2k}{k-q} - \left [\binom{2k}{k+q-1-2m} +\binom{2k}{k-q-1-2m}\right ] \cos (2\theta) & 2m+1+q \le k< 4m+2-q \nonumber \\
... & ...
\end{cases}
& \nonumber \\
&\tfrac{1}{2m+1}{\sum}_{\ell=0}^{2m} [\sin \tfrac{\theta+2\pi\ell}{2m+1}]^{2k+1} \sin (2q \tfrac{\theta+2\pi\ell}{2m+1}) =\\
&=\tfrac{(-)^{q+m+1}}{2^{2k+1}} 
\begin{cases} 0 & k<m-q \nonumber\\
\binom{2k+1}{k-m+q} \cos\theta & m-q\le k< m+q \nonumber \\
\left [ \binom{2k+1}{k-m+q} - \binom{2k+1}{k-m-q}\right ]\cos \theta & m+q\le k< 3m+1-q \nonumber\\
... & ... \nonumber
\end{cases} 
\end{align}
{\bf Case $\mathbf{n=2m+1}$ and $\mathbf{p=2q+1}$}:
\begin{align*}
\tfrac{1}{2m+1}{\sum}_{\ell=0}^{2m}& \cos[x\sin \tfrac{\theta+2\pi\ell}{2m+1}] \cos[ (2q+1) \tfrac{\theta+2\pi\ell}{2m+1}] = 
[J_{2m-2q}(x)+J_{2m+2q+2}(x)]\cos\theta+...\\
-\tfrac{1}{2m+1}{\sum}_{\ell=0}^{2m}&\sin [x \sin\tfrac{\theta+2\pi\ell}{2m+1}] \sin [(2q+1) \tfrac{\theta+2\pi\ell}{2m+1}] = J_{2q+1}(x)+\\
&+[-J_{4m-2q+1}(x)+J_{4m+2q+3}(x)] \cos(2\theta)+ ...
\end{align*}
\begin{align}
&\tfrac{1}{2m+1}{\sum}_{\ell=0}^{2m} \sin [\tfrac{\theta+2\pi\ell}{2m+1}]^{2k}  \cos[ (2q+1) \tfrac{\theta+2\pi\ell}{2m+1}] =\\
&=\tfrac{(-)^{m+q}}{4^k}\begin{cases} 0 & k<m-q  \nonumber\\
\binom{2k}{k-m+q} \cos\theta & m-q\le k < m+q+1 \nonumber \\
\left[ \binom{2k}{k-m+q} -\binom{2k}{k-m-q-1} \right ] \cos\theta & m+q+1\le k < 3m+q+2 \nonumber\\
... & ...
\end{cases}
\end{align}
\begin{align}
&\tfrac{1}{2m+1}{\sum}_{\ell=0}^{2m} \sin [\tfrac{\theta+2\pi\ell}{2m+1}]^{2k+1}  \sin[ (2q+1) \tfrac{\theta+2\pi\ell}{2m+1}] =\\
&=\tfrac{(-)^q}{2^{2k+1}} \begin{cases} 0 & k<q  \nonumber\\
\binom{2k+1}{k-q} & q\le k <  2m-q    \nonumber\\
\binom{2k+1}{k-q} - \binom{2k+1}{k+q-2m} \cos(2\theta) & 2m- q\le k <  2m+q+1    \nonumber\\
\binom{2k+1}{k-q} - \left[ \binom{2k+1}{k+q-2m} + \binom{2k+1}{k-q-1-2m}\right ]\cos(2\theta) & 2m+ q+1\le k <  4m-q+1    \nonumber\\
... & ... \\
\end{cases}
\end{align}

 \subsection*{Data availability}
Data sharing is not applicable to this article as no new data were created or analyzed in this study.

\end{document}